\documentclass[12pt]{amsart}
\usepackage[T2A]{fontenc}
\usepackage[cp1251]{inputenc}
\usepackage[english]{babel}
\usepackage{amsmath,latexsym,amsthm,amsfonts,tipa,upgreek}
\usepackage{amssymb,amscd,graphpap, stmaryrd}

\begin{document}

\title[LUR - нормы на компактах Федорчука]{Существование LUR - норм на пространствах \\ непрерывных 
функций на компактах Федорчука}
\author{С.П.Гулько, А.В.Иванов, С.Троянский, М.С.Шуликина}

\thanks{The study was carried out with the financial support of the Russian Foundation for Basic
Research in the framework of the scientic project N 17-51-18051. 
The third author was supported by Bulgarian National Scientific Fund under Grant DEMI/Russia, 01/06/23.06.2017}

\maketitle

\begin{abstract}
В данной статье построена LUR-норма на пространстве $C(X)$, которая является полунепрерывной снизу относительно топологии поточечкой сходимости для любого компакта Федорчука $X$ высоты 3.
\end{abstract}


 В данной заметке рассматриваются примеры компактов Федорчука $X$, для которых банахово пространство $C(X)$ имеет эквивалентную $LUR$-норму, которая является полунепрерывной снизу относительно топологии поточечной сходимости. Напомним, что компакт $X$ является компактом Федорчука, если его можно получить как предел обратного спектра, который начинается с одноточечного пространства и у которого все соседние проекции являются вполне замкнутыми отображениями с метризуемыми слоями. Спектральной высотой компакта называется наименьшая длина обратного спектра, в который этот компакт можно разложить. В этой заметке нас будут интересовать компакты со спектральной высотой 3. 
Из определении компактов Федорчука следует , что компакты
спектральной высоты 1 ---  это точки, метрические компакты
имеют высоту 2  и неметризуемые компакты имеют
спектральную высоту не меньше 3. На самом деле компакт
Федурчука $Х$ высоты 3 определяется следующим
образом --- это неметризуемые компакты, которые допускают вполне замкнутое отображение на метрический компакт и имеющие метрические слои.
Мир компактов Федорчука высоты 3 достаточно обширен. Им фактически посвящена обзорная статья С.Ватсона [5]. Кроме того, ряд типичных компактов Федорчука имеют высоту 3. Например, таков знаменитый пример компакта, размещенный  в книге [4] на стр.321. 

Основным результатом данной статьи является следующая теорема

\textbf{Теорема} Если $X$ является компактом Федорчука спектральной высоты 3, то для пространства $C(X)$ всех непрерывных вещественных функций на $X$ существует эквивалентная локально равномерно выпуклая норма (LUR-норма), которая полунепрерывна снизу относительно топологии поточечной сходимости. 

Нам будет удобнее всего использовать следующее определение вполне замкнутого отображения.

\textbf{Предложение 1} [2]. 
Отображение $f:X\to Y$ является вполне замкнутым тогда и только тогда, когда пересечение $f(F_1)\cap f(F_2)$ конечно для любых замкнутых непересекающихся множеств $F_1$ и $F_2$ в $X$.

Другую характеризацию вполне замкнутых отображений дает следующее утверждение.

{\bf Предложение 2.} {\it Отображение компактов $f:X\to Y$ вполне замкнуто тогда и только тогда, когда для любого непрерывного отображения $g:X\to K$ в метрический компакт $K$ и любого $\varepsilon>0$ множество $H_{g,\varepsilon}=\{y\in Y: diam(g(f^{-1}(y)))\geq\varepsilon\}$ конечно.}

{\bf Доказательство.} Необходимость. Пусть $f$ вполне замкнуто. Это значит, что для любых непересекающихся замкнутых подмножеств $A,B\subset X$ пересечение $f(A)\cap f(B)$ конечно. Предположим, что существует непрерывное отображение $g:X\to K$ в метрический компакт $K$ и $\varepsilon>0$  такие, что $|H_{g,\varepsilon}|\geq\omega$. Рассмотрим бесконечное индексированное элементами $y$ семейство замкнутых подмножеств $K$:
$$
G=\{g(f^{-1}(y)):y\in H_{g,\varepsilon}\}.
$$
Пусть $F$ -- точка полного накопления семейства $G$ в пространстве $\exp(K)$ непустых замкнутых подмножеств $K$ с метрикой  Хаусдорфа. Очевидно, что $diam(F)\geq\varepsilon$. Возьмем в $F$ две различные точки $x_1,x_2$ и их окрестности $Ox_1,Ox_2$ с непересекающимися замыканиями. В силу выбора $F$ существует бесконечное подмножество $D\subset H_{g,\varepsilon}$ такое, что $g(f^{-1}(y))\cap Ox_i\not=\emptyset$ для любого $y\in D$ при $i=1,2$. Рассмотрим непересекающиеся замкнутые подмножества $A_i=g^{-1}[Ox_i],\ i=1,2$ в пространстве $X$. Имеем $D\subset f(A_1)\cap f(A_2)$, что противоречит вполне замкнутости $f$.

Достаточность. Пусть $f$ не является вполне замкнутым, то есть в $X$ существуют замкнутые непересекающиеся подмножества $A,B$, для которых пересечение $E=f(A)\cap f(B)$ бесконечно. Пусть $g:X\to [0,1]$ -- непрерывная функция на $X$, разделяющая $A$ и $B$. Тогда множество $H_{g,1}$ содержит в себе $E$ и, следовательно, бесконечно. $\Box$

Далее нам понадобятся следующие определения из [1].

{\bf Определение.} {\it Для любого вполне замкнутого отображения $f:X\to Y$ отображение послойной осцилляции $\Omega_f: C(X)\to c_0(Y)$ определяется по формуле:
$$
\Omega_f(g)(y)=diam(g(f^{-1}(y))),
$$
где $g\in C(X)$.}

В силу предложения 2 определение величины $\Omega_f$ является корректным.

Для доказательства основной теоремы статьи нам понадобятся еще два утверждения. 

{\bf Предложение 3.} {\it Отображение послойной осцилляции $\Omega_f: C(X)\to c_0(Y)$ ко-$\sigma$-непрерывно для любого вполне замкнутого отображения $f: X\to Y$ на метрический компакт $Y$ с метризуемыми слоями $f^{-1}(y), y\in Y$.}

Термин ко-$\sigma$-непрерывно в данном случае означает, что для каждого $g\in C(X)$ существует сепарабельное подмножество $Z_g$ пространства $C(X)$ такое, что $g$ лежит в замыкании (по норме) линейной оболочки объединения $Z_{g_n}$ как только $\lim\Omega_f(g_n)=\Omega_f(g)$.

{\bf Доказательство.} Воспользуемся критерием ко-$\sigma$-непрерывности из [1] (теорема 1.7). Пусть $g\in C(X)$. Положим
$$
Z_g=\{h\in C(X): supp\ \Omega h\subset supp\ \Omega g\}.
$$
Множества $Z_g$ сепарабельны для любого $g$. 

Пусть $\lim\Omega_f(g_n)=\Omega_f(g)$. Покажем, что $g\in [\cup Z_{g_n}]$.

Рассмотрим разбиение $X$, нетривиальными элементами которого являются слои $f^{-1}(y)$ при $y\not\in supp\ \Omega g$. Фактор-пространство по этому разбиению обозначим через $X_1$. Из результатов [2] следует, что $X_1$ хаусдорфово.  Пусть $h_1: X\to X_1$ -- факторное отображение и отображение $f_1:X_1\to Y$ таково, что $f_1\circ h_1=f$. Тогда определено непрерывное отображение $g'=g\circ h_1^{-1}: X_1\to R$ и $g=g'\circ h_1$.

Фиксируем $\varepsilon>0$, и пусть $K=\{y: diam(g(f^{-1}(y)))\geq\varepsilon/2\}$, $|K|<\omega$. Рассмотрим разбиение $X_1$, нетривиальными элементами которого являются множества $f_1^{-1}(y)$ при $y\in supp\ \Omega g\setminus K$. Фактор-пространство по этому разбиению обозначим через $X_2$, а факторное отображение $X_1$ на $X_2$ -- через $h_2$. Пусть $f_2: X_2\to Y$ таково, что $f_2\circ h_2=f_1$. Для каждого $y\in supp\ \Omega g\setminus K$ выполнено $diam (g'(f_1^{-1}(y)))<\varepsilon/2$. Пусть $O_y$ -- окрестность множества $f_1^{-1}(y)$ в $X_1$ такая, что
$$
diam(g'(O_y))<2\ diam(g'(f_1^{-1}(y))).
$$
Занумеруем точки $supp\ \Omega g\setminus K$ натуральными числами: $supp\ \Omega g\setminus K=\{y_n:n\in N\}$ и обозначим через $y_n'$ единственную точку $X_2$, для которой $f_2(y_n')=y_n$.

Зафиксируем метрику $d$ на $X_2$, совместимую с топологией ($X_2$ метризуемо) и построим по рекурсии открытые множества $O_n,\ n\in N$, в $X_2$ следующим образом.

Шаг 1. По построению множеств $O_y$, $y_1'\in h^\sharp_2O_{y_1}$. Множество $f_2^{-1}(supp\ \Omega g\setminus K)$ счетно, следовательно, нульмерно. Поэтому существует окрестность $O_1$ точки $y_1'$ такая, что $[O_1]\subset h_2^\sharp O_{y_1}$, $fr\ O_1\cap(f_2^{-1}(supp\ \Omega g\setminus K))=\emptyset$ и $diam[O_1]<1$.

Шаг $n$. Если $y_n'\in \cup_{k<n}O_k$, то положим $O_n=\emptyset$. В противном случае выберем $O_n$ так, что  $[O_n]\subset h_2^\sharp O_{y_n}$, $fr\ O_n\cap(f_2^{-1}(supp\ \Omega g\setminus K))=\emptyset$, $diam[O_n]<1/n$ и $[O_n]\cap(\cup_{k<n}[O_k])=\emptyset$.

В результате построения мы получим семейство $O_n,\ n\in N$ такое, что
$$
\cup O_n\supset f_2^{-1}(supp\ \Omega g\setminus K).
$$

Пусть $U=h_2^{-1}(\cup O_n)$ и $F=X_1\setminus U$. Отображение $h_2$ взаимно однозначно на $F$ и $h_2(F)=G=X_2\setminus\cup O_n$. Положим
$$
r=g'\circ(h_2^{-1}|_G):G\to R.
$$
По построению $diam(g'(O_{y_n}))<2\ diam(g'(f_1^{-1}(y_n)))$ для каждого $n$. Положим $a_n=\inf(g'(O_{y_n})),\ b_n=\sup(g'(O_{y_n}))$, $b_n-a_n<2\ diam(g'(f_1^{-1}(y_n)))<\varepsilon$. Таким образом, $r(fr(O_n))\subset[a_n,b_n]$. Продолжим отображение $r|_{fr( O_n)}$ до непреывного отображения $r_n:[O_n]\to[a_n,b_n]$.

Определим теперь отображение $p':X_2\to R$ следующим образом: при $x\in G$ $p'(x)=r(x)$,  при $x\in O_n$ $p'(x)=r_n(x)$.

Покажем, что отображение $p'$ непрерывно. Если $x\in\cup O_n$, то непреывность $p'$ в точке $x$ очевидна.

Пусть $x\in G$ и $x\not\in\cup[O_n]$. Фиксируем $a>0$. Пусть $z$ -- единственная точка $X_1$, для которой $h_2(z)=x$. В силу непрерывности $g'$, существует окрестность $Oz$ такая, что для любого $z'\in Oz$ $|g'(z)-g'(z')|<a/5$.   Имеем $x\in h_2^\sharp(Oz)$. Пусть число $k\in N$ таково, что $O(x,1/k,d)\subset h_2^\sharp(Oz)$. Положим $Ox=O(x,1/2k,d)\setminus\cup_{n<2k}[On]$. 

Покажем, что для любого $x'\in Ox$ $|p'(x)-p'(x')|<a$. Если $x'\in G$, то имеется единственная точка $z'\in Oz$ такая, что $h_2(z')=x'$, и $p'(x')=g'(z')$. Следовательно, $|p'(x)-p'(x')|=|g'(z)-g'(z')|<a/5$.

Остается случай, когда $x'\in O_n$ где $n\geq 2k$. Поскольку $diam(O_n)<1/n$, $O_n\subset O(x,1/k,d)\subset h_2^\sharp(Oz)$. Следовательно, $f_1^{-1}(y_n)\subset Oz$. Значит, $diam(g'(f_1^{-1}(y_n)))<2a/5$. По построению множеств $O_n$ для любой точки $z'\in h_2^{-1}(x')$ $|g'(z')-p'(x')|<4a/5$. Но $g'(z)=p'(x)$, откуда следует, что $|p'(x)-p'(x')|<a$ -- непрерывность $p'$ в точке $x$ доказана.

Осталось рассмотреть случай, когда $x\in fr(O_n)$ для некоторого $n$. В этом случае представим $X_2$ в виде объединения двух замкнутых множеств:
 $$
X_2=[O_n]\cup (X_2\setminus O_n).
$$
Ограничение $p'$ на $[O_n]$ непрерывно в точке $x$ по определению $p'$. Непрерывность ограничения $p'$ на $(X_2\setminus O_n)$ в точке $x$ доказывается дословным повторением проведенных выше рассуждений с заменой $X_2$ на $(X_2\setminus O_n)$, $X_1$ -- на $X_1\setminus h_2^{-1}(O_n)$, а $O_n$ -- на пустое множество. Из непрерывности $p'$ на $[O_n]$  и $(X_2\setminus O_n)$ следует непрерывность на $X_2$.

Положим теперь $p=h_1\circ h_2\circ p': X\to R$. По построению имеем

1)$||p-g||\leq\varepsilon$;

2) $supp \ \Omega p=K\subset supp\ \Omega g$.

Из $\lim\ \Omega g_n=\Omega g$  следует, что существует $n_0$ такое, что $K\subset supp\ \Omega g_n$ при $n>n_0$. Значит, $p\in \cup Z_{g_n}$, что и требовалось. $\Box$

{\bf Предложение 4.} {\it Отображение послойной осцилляции $\Omega_f$ является $\sigma$-срезочно (slicely) непрерывно.}

{\bf Доказательство.} В [1] доказано (Следствие 1.21), что всякое локально ограниченное отображение $\Phi$ из нормированного пространства $X$ в $c_0(\Gamma)$ такое, что для любого $\gamma\in\Gamma$ функция $\delta_\gamma\circ \Phi:X\to R$ неотрицательна и выпукла, является $\sigma$-срезочно непрерывным. Выполнение всех указанных условий для отображения $\Omega_f:C(X)\to c_0(Y)$ устанавливается непосредственной проверкой.

Отображение $\Omega_f$ локально ограничено относительно топологии на $C(X)$, порожденной нормой.

Функция $\delta_y\circ\Omega_f$ определяется равенством:
$$
\delta_y\circ\Omega_f(g)=diam(g(f^{-1}(y))).
$$
Таким образом, эта функция неотрицательна и выпукла, поскольку
$$
\delta_y\circ\Omega_f(tg_1+(1-t)g_2)\leq\delta_y\circ\Omega_f(tg_1)+\delta_y\circ\Omega_f((1-t)g_2)=
$$
$$
=t(\delta_y\circ\Omega_f(g_1))+(1-t)(\delta_y\circ\Omega_f(g_2))
$$
при $g_1,g_2\in C(X)$ и $t\in (0,1)$. $\Box$

Основная теорема данной статьи теперь следует из перечисленнных предложений, ср. [1].

\end{document}